\documentclass[12pt]{amsart}
\usepackage{amsfonts}
\usepackage{amssymb,amsthm,amscd,amstext,amsmath,mathrsfs,latexsym}
\newtheorem{theorem}{Theorem}

\newcommand{\dai}{|\Delta|a}

\begin{document}

\title[]
{The FOUR NUMBER GAME}

\author{BY BENEDICT FREEDMAN}
\email{}
\address{}

\thanks{}

\maketitle

%\Large

THE mathematical curiosity called The \lq\lq Four Number Game" was shown to me by Prof. Jekuthiel Ginsburg in 1938 with instructions to work out the theory of it as an exercise in algebra.  I hope he will excuse me for handing in my homework ten years late.

The game is played as follows.  Take any four positive integers, e.g., $8,17,3,107$.  Under each  number write the numerical difference between it and the following number.  Thus:

$$\begin{matrix}
8&17&3&107\\
9&14&104&99
\end{matrix}$$

Two points should be noted.  The positive value of the difference is always taken, $9=+|8-17|$, $14=+|17-3|$, etc.  To obtain the last difference the first number is considered to follow the fourth.  Therefore $99=+|107-8|$.

Now we have a new set of four numbers which may be differenced in the same way.  Repeating, we get the sequence:

$$\begin{matrix}
A_0=&8&17&3  &107\\
A_1=&9&14&104&99\\
A_2=&5&90&5&90\\
A_3=&85&85&85&85\\
A_4=&0&0&0&0
\end{matrix}$$

When we reach a set consisting of four zeros, the game is at an end.  Experiment reveals that any set of four numbers leads to four zeros; and it turns out to be next to impossible to pick a set that will last for more than six or eight differencings.  Even as widely spaced a quartet as $1, 11, 130, 1760,$ holds out for only six steps.\footnote{$$\begin{matrix}
A_0=&1&11&130  &1760\\
A_1=&10&119&1630&1759\\
A_2=&109&1511&129&1749\\
A_3=&1402&1382&1620&1640\\
A_4=&20&238& 20& 238\\
A_5=&218&218&218&218\\
A_6=&0&0&0&0
\end{matrix}$$}
The reader is invited to better this by trial.  When you have just about come to the conclusion that eight differencings is the maximum, it is infuriating to run across a set like $0, 6, 17, 37$ which survives twelve steps.
$$\begin{matrix}
A_0=&0&6&17  &37\\
A_1=&6&11&20&37\\
A_2=&5&9&17&31\\
A_3=&4&8&14&26\\
A_4=&4&6&12&22\\
A_5=&2&6 &10 &18\\
A_6=&4&4&8&16\\
A_7=&0&4&8&12\\
A_8=&4&4&4&12\\
A_9=&0&0&8&8\\
A_{10}=&0&8&0&8\\
A_{11}=&8&8&8&8\\
A_{12}=&0&0&0&0
\end{matrix}$$

It is natural to ask two questions:

\begin{itemize}
\item[1.] Given \emph{any} four numbers as the original set, will we reach the set of four zeros after a finite number of differencings?
\item[2.] If so, is there a maximum number of differencings, or can sets be found that will outlast $n$ steps where $n$ is as large as you please?
\end{itemize}

Although for convenience the four numbers of a set are written on a line, writing them around a circle will exhibit more clearly certain cyclic properties.  For example, the sets $(a,b,c,d), (b,c,d,a), (c,d,a,b), (d,a,b,c)$ are for all purposes equivalent; the only difference being in the member chosen as initial.  The differences of the second set are the corresponding differences of the first set shifted one space to the left; and so on.  Also the sets obtained by reading any of the above backwards are equivalent to them.

Suppose that the original set consists of three even numbers followed by an odd:
$A_0=e\;e\;e\;o$.  Without any further information about the magnitude of the numbers, the even- or odd-ness of the successive differences can be ascertained, since the difference of two even numbers is even, of two odds is even, and of even and odd is odd.
$$\begin{matrix}
A_0=&e&e&e&o\\
A_1=&e&e&o&o\\
A_2=&e&o&e&o\\
A_3=&o&o&o&o\\
A_4=&e&e&e&e
\end{matrix}$$

Notice that the fourth set of differences is composed solely of even numbers.  This is true no matter what the original set is.  Because of the cyclic equivalences discussed above, only $(eeee), (eeeo), (eeoo), (eoeo), (eooo),$ and $(oooo)$ need be tested.  Any other combination is equivalent to one of these.  Therefore we have the result that for any set of four numbers, the fourth set of differences will be all even.

Since the members of $A_4$ are all even, we are led to consider the relation between $A_4$ and $\frac{1}{2}A_4$, the set whose members are halves of those of $A_4$.  This relationship is preserved throughout differencing.  That is, the $r$th differences of $A_4$ are double the corresponding $r$th differences of $\frac{1}{2}A_4$.  But we know that the fourth differences of $\frac{1}{2}A_4$ are even; therefore the fourth differences of $A_4$, i.e., $A_8$, being double these, are divisible by four.  In like manner, as the old algebras say, $A_{12}$ is divisible by eight, $A_{16}$ by sixteen, and in general $A_{4h}$ is divisible by $2^h$.

Consider any original set $A_0=a\;b\;c\;d\;$.  Suppose that none of the numbers is larger than $1000$.  Our method of numerical differencing can never produce a larger number than the largest in the previous set.  Therefore the members of the fortieth difference of $A_0$ are none of them greater than $1000$.  But we have just shown that $A_{40}$ is divisible by $2^{10}$, i.e., by $1024$.  The only number leas than $1000$ and divisible by $1024$ is zero.  Therefore $A_{40}$ is all zeros.  It is not hard to generalize the proof with any number $n$ taking the place of $1000$.  Then we can answer \lq\lq Yes" to question one.

The results for sets of four numbers are of course special cases of theorems for sets of $k$ numbers which we will now develop.  Using the operator $|\Delta|a_i=|a_{i+1}-a_i|$, and adopting the convention that subscripts are to be reduced modulo $k$ (e.g. $|\Delta|a_k=|a_{k+1}-a_k|=|a_1-a_k|)$, we can represent the derived sets as follows:
$$\begin{matrix}
A_0=&a_1&a_2&a_3&...&a_k\\
A_1=&\dai_1&\dai_2&\dai_3&...&\dai_k\\
A_2=&|\Delta|^2a_1&|\Delta|^2a_2&|\Delta|^2a_3&...&|\Delta|^2a_k\\
\vdots &&&&&\\
A_n=&|\Delta|^n a_i &i=1,2,...k,&&&
\end{matrix}$$
where $|\Delta|^n a_i=|[ |\Delta|^{n-1} a_{i+1}-|\Delta|^{n-1} a_i] |$.

The relation between $|\Delta|$ and the ordinary difference operator, $\Delta$, is simple:
\begin{equation}
\dai_i=|\Delta a_i|
\end{equation}
since each side $=|a_{i+1}-a_i|$.  Also the difference $[|a_{i+1}-a_i|-(a_{i+1}-a_i)]$ is always even.  Therefore,
$$|\Delta|a_i\equiv \Delta a_i\; \textrm{mod} \; 2$$
By Mathematical Induction, using the definition of $|\Delta|^n$, we get
\begin{equation}
|\Delta|^na_i=\Delta^n a_i \;\textrm{mod} \; 2
\end{equation}

Expanding the right side by Newton's Advancing Difference Formula,
\begin{equation}
|\Delta|^na_i=[a_{i+n}-{n \choose 1}a_{i+n-1}+{n \choose 2}a_{i+n-2}-...+(-1)^n a_i]\; \textrm{mod} \; 2
\end{equation}
it being understood, as before, that the subscripts are to be reduced modulo $k$.  Formula $(3)$ may be written symbolically more compactly
\begin{equation}
|\Delta|^n a_i\equiv (a-1)^n a_i \; \textrm{mod} \; 2.
\end{equation}

Take $n=2^r$, and using the fact that ${2^r \choose j}$ is always even for $0< j< 2^r$,\footnote{\it{Proof:}
$$\begin{matrix}
{2^r \choose j}=&\prod_{\alpha=0}^{j-1} \frac{2^r-\alpha}{\alpha+1} &0<j< 2^r\\
&=\frac{2^r}{j}\prod_{\alpha=1}^{j-1}\frac{2^r-\alpha}{\alpha } &\textrm{(rearranging the factors)}
\end{matrix}$$
Every factor $\frac{(2^r-\alpha)}{\alpha}$ can be reduced to a fraction with an odd denominator.  Since the highest power of $2$ in $j$ must be less than $r$, the product ${2^r \choose j}$ is even.  That is, ${2^r \choose j}\equiv 0 \; \textrm{mod} \; 2$ for $0< j < 2^r$.
}
we arrive at a useful result:
\begin{equation}
|\Delta|^{2^r}a_i\equiv (a_{i+2^r}+a_i) \; \textrm{mod} \; 2.
\end{equation}
If $k$, the number of members of the set is a power of $2$, $k=2^r$, $(5)$ becomes
\begin{equation}
|\Delta|^k a_i\equiv (a_{i+k}+a_i)\equiv (a_i+a_i)\equiv 0 \; \textrm{mod} \; 2, k=2^r.
\end{equation}
That is, if $A_0$ has $k=2^r$ members, every member of $A_k$ is even.  A check of our previous examples will show that in every case $A_4$ is composed of even numbers.  Here, $k=4=2^2$.

Before proceeding, we need the following lemma:

If $a_i=\lambda b_i$, where $\lambda$ is a positive integer, then
\begin{equation}
|\Delta|^n a_i=\lambda |\Delta|^n b_i.
\end{equation}

Proof is by mathematical induction, starting from:
\begin{center}
$a_i=\lambda b_i$ implies $|a_{i+1}-a_i|=\lambda |b_{i+1}-b_i|.$
\end{center}
In words, if each member of one set is a fixed multiple of the corresponding member of another set, then the same relation holds between their respective differences.\footnote{Example:$$\begin{matrix}
[3&12&6&6]&=3&\times &[1&4&2&2]\\
[9&6&0&3]&=3&\times &[3&2&0&1]\\
[3&6&3&6]&=3&\times &[1&2&1&2]\end{matrix}$$
}

Returning to the main thread of our development, if $k=2^r$, then by $(6)$
\begin{center}
$|\Delta|^ka_i\equiv 0 \; \textrm{mod} \; 2 \;\;\; i=1,2,...k.$
\end{center}
Set $b_i=\frac{1}{2} |\Delta|^k a_i.$  Since $|\Delta|^ka_i$ is even, $b_i$ are integers.  Then
\begin{center}
$|\Delta|^ka_i=2b_i \;\;\; i=1,2,...k.$
\end{center}
By $(7)$, $|\Delta|^{2k}a_i=|\Delta|^k(2b_i)=2|\Delta|^kb_i.$  But, by $(6)$, since $k=2^r$
\begin{center}
$|\Delta|^kb_i\equiv 0 \; \textrm{mod} \; 2$
\end{center}
\begin{center}
$\therefore |\Delta|^{2k} a_i\equiv 0 \; \textrm{mod} \; 4.$
\end{center}

Repeating this process, we conclude
\begin{equation}
|\Delta|^{hk}a_i\equiv 0 \; \textrm{mod} \; 2^k \;\;\; i=1,2,...k, k=2^r.
\end{equation}

Now, let $c$ be a number greater than every member of the original set, $c> a_i$.  Since the process of differencing involves subtraction only, $c>|\Delta|^{hk} a_i$.  Take $h$ such that $2^h>c$.  Then $|\Delta|^{hk}a_i < 2^h$ and $|\Delta|^{hk}a_i\equiv 0 \; \textrm{mod} \; 2^h$ which means that $|\Delta|^{hk}a_i=0,$ i.e., every member of $A_{hk}$ equals zero.  Thus, we have

\begin{theorem}
If the number of members of the original set is a power of $2$ $(k=2^r)$, then the process of differencing leads ultimately to a set where every member is zero.
\end{theorem}

What if the number of the members of the original set is \emph{not} a power of $2$?  We shall show that only in special cases does differencing lead to an all-zero set.

Let $A_0=a_1\;a_2\;a_3... a_k$, and let $k=2^{\alpha} K$, where $\alpha\geq 0$ and $K$ is an odd number greater than $1$.  Suppose that the $n$th derived set, $A_n=|\Delta|^na_i$, is all zeros.  We shall prove that this imposes a very restrictive condition on $a_i$.

By a familiar theorem in the Theory of Numbers,
\begin{center}
$2^{\phi(K)}\equiv 1\; \textrm{mod} \; K$, since $2$ is prime to $K$.
\end{center}
Raising both sides to an arbitrary power $\lambda$,
\begin{center}
$2^{\lambda \phi(K)}\equiv 1\; \textrm{mod} \; K.$
\end{center}
Multiplying the congruence by $2^\alpha$,
\begin{center}
$2^{\lambda \phi(K)+\alpha}\equiv 2^\alpha\; \textrm{mod} \; (2^\alpha K).$
\end{center}
Or

\begin{equation}
2^{\lambda \phi(K)+\alpha}\equiv 2^\alpha\; \textrm{mod} \; k.
\end{equation}

Since $\lambda$ is arbitrary, we can take it large enough so that
\begin{center}
$2^{\lambda \phi(K)+\alpha}\geq n.$
\end{center}
Then, since $|\Delta|^na_i=0$ by hypothesis,
\begin{equation}
|\Delta|^ma_i=0, \textrm{where}\; m \; \textrm{stands for}\; 2^{\lambda \phi(K)+\alpha}.
\end{equation}

But $m$ is a power of $2$.  Therefore, by $(5)$,
\begin{center}
$|\Delta|^m a_i\equiv a_{i+m}+a_i \; \textrm{mod} \; 2.$
\end{center}
Combining with $(10)$,
\begin{center}
$a_{i+m}+a_i\equiv 0 \; \textrm{mod} \; 2,$
\end{center}
But subscripts are to be reduced modulo $k$, and by $(9)$, $m\equiv 2^\alpha \; \textrm{mod} \; k$.
\begin{equation}
a_{i+2^\alpha}+a_i\equiv 0 \; \textrm{mod} \; 2
\end{equation}

or
\begin{center}
$a_{i+2^\alpha}\equiv a_i \; \textrm{mod} \; 2$ for $i=1,2,...k.$
\end{center}
This restrictive condition is obviously not satisfied by most original sets $A_0$.  We can state a second theorem:
\begin{theorem}
If the number of members of the original set is \emph{not} a power of $2$, the process of differencing does not in general lead to a set of all zeros.
\end{theorem}

An example will help make this clear,
$$\begin{matrix}
A_0&=&2&5&9\\
A_1&=&3&4&7\\
A_2&=&1&3&4\\
A_3&=&2&1&3\\
A_4&=&1&2&1\\
A_5&=&1&1&0\\
A_6&=&0&1&1\\
A_7&=&1&0&1\\
A_8&=&1&1&0\\
&&etc.&&
\end{matrix}$$
Note that $A_8$ is identical with $A_5$, and that $A_5, A_6, A_7$ repeat from then on.  This is true in general.  For if $A_0$ does not lead to an
all-zero set, it must lead to an infinite sequence of sets.  But as subtraction is the only operation employed, no number greater than the greatest in $A_0$ can be introduced.  There are only a finite number of arrangements of numbers less than a given number (taken $k$ at a time).  Therefore, at some stage repetition must take place.

We know now that an original set of four numbers leads to four zeros, and casual trial indicates that the majority of such sets arrive at four zeros in less than ten differecings.  However, we shall show that it is possible to construct sets that will outlast any assigned number of differencings.

Represent by $\{a_i\}$ the number of steps it takes to bring the set $a_i$ to four zeros.  That is, $n=\{a_i\}$ is the least integer for which
\begin{center}
$|\Delta|^na_i=0 \;\;\; i=1,2,3,4.$
\end{center}
As a corollary of formula $(7)$, we have
\begin{equation}
\{\lambda a_i\}=\{a_i\},
\end{equation}
where $\lambda$ is a positive integer.

Since $|\Delta|(a_i+\delta)=|(a_{i+1}+\delta)-(a_i+\delta)|=|\Delta|a_i,$
\begin{equation}
\{a_i+\delta\}=\{a_i\} \;\;\;  \delta\geq 0
\end{equation}
Combining $(12)$ and $(13)$, we get
\begin{equation}
\{\lambda a_i+\delta\}=\{a_i\}\;\;\; \lambda, \delta \geq 0.
\end{equation}

If we multiply each member of a set by a constant and add a constant to each, the new set has the same \lq\lq life expectancy," so to speak, as the original.

Let us now investigate the condition under which one set can be the differences of another.  Suppose
\begin{equation}
b_i=\dai_i \;\;\; i=1,2,...k.
\end{equation}
Consider the numbers $B_i=\Delta a_i$.  Clearly,
$$\sum_{1}^{k}B_i=\sum_{1}^{k}\Delta a_i=a_{k+1}-a_1=a_1-a_1=0.$$
Separating the positive $B_i$ from the negative, $\sum B_i=\sum-B_j$, where $B_i\geq 0$ and $B_j<0$.  By $(1)$, $b_i=|B_i|$.  Therefore,
\begin{equation}
\sum b_i=\sum b_j
\end{equation}
that is, for the $b$'s to be the differences of another set, the sum of some of the $b$'s must equal the sum of the remainder. \footnote{For example $(2,3,7,7)$ can \emph{not} be a set of differences, since there is no way to group them into two equal sums.  On the other hand $(2,3,7,6)$ \emph{is} a set of differences since $2+7=3+6$.  It is the difference of $(5,7,4,11)$ to name only one.  $(2,3,7,12)$ also satisfies the criterion $(2+3+7=12)$ and is the differences of $(0,2,5,12)$.}

We are ready to answer question two, first for $k=4$, and then in general.  Consider a set $(a_1\; a_2\; a_3\; a_4)$, for which
\begin{itemize}
\item[$(\alpha)$] $a_1=0$; $a_2\; a_3\; a_4$ are not all $0$
\item[$(\beta)$] $a_3\geq a_2\geq 0$
\item[$(\gamma)$] $a_4\geq a_3+a_2.$
\end{itemize}

We will find a set which will last one more step than $a_i$.  First, if $a_4=a_3+a_2+a_1$, condition $(16)$ is satisfied and the $a'$s are the differences of another set.  If not, using $(14)$ we can construct a set with the same \lq\lq life" as $a_i$ which does satisfy $(16)$.  To do this, place
$$a_i'=\lambda (a_i) +\delta \;\;\; i=1,2,3,4$$
$$a_4'=a_3'+a_2'+a_1'.$$
Eliminating the primes,
$$\lambda a_4+\delta=\lambda a_3+\delta+\lambda a_2+\delta +\lambda a_1 +\delta.$$
Solving for $\delta$
\begin{equation}
\delta=\frac{\lambda}{2}(a_4-a_3-a_2-a_1).
\end{equation}
If we take $\lambda=2$, $\delta$ will be a positive integer since $a_i$ satisfy restrictions $\alpha, \beta,\gamma.$  Therefore $a_i'$ are positive integers and also (since $a_4'=a_3'+a_2'+a_1'$) the differences of a new set, $b_i$, for which
$$b_1=0$$
$$b_2=a_1'$$
$$b_3=a_1'+a_2'$$
$$b_4=a_1'+a_2'+a_3' (=a_4').$$
It is easy to check that $a_i'=|\Delta|b_i.$  Furthermore, $b_i$ also satisfy restrictions $\alpha, \beta, \gamma.$\footnote{$\alpha$ and $\beta$ can be proved by inspection.  For $\gamma$, note that $b_4\geq b_3+b_2$ reduces to $$a_1'+a_2'+a_3'\geq a_1'+a_2'+a_1'$$ $$a_3'\geq a_1'$$  $$\lambda a_3+\delta\geq \lambda a_1 +\delta \;(\textrm{eliminating primes)}$$ $$a_3\geq a_1=0 \;(\textrm{which is given}).$$}

By $(14)$,
$$ \{a_i'\}=\{a_i\}.$$
And since $a_i'$ are the differences of $b_i$,
$$ \{b_i\}=\{a_i'\}+1.$$
Or,
\begin{equation}
\{b_i\}=\{a_i\}+1\;\;\; i=1,2,3,4.
\end{equation}
Thus, if a set $a_i$ satisfies the restrictions $\alpha, \beta, \gamma,$ we can construct a set which will last one step longer.  But as we have seen, this new set will also satisfy the restrictions $\alpha, \beta, \gamma$.  Therefore the process may be repeated on it, and so on indefinitely.  That means it is possible to construct sets of four which will outlast any specified number of differencings.

\emph{Example:}

$$\begin{matrix}
a_i=&0&0&0&1&   &\{a_i\}=4\\
a_i'=&1&1&1&3&(\lambda=2;\delta=1)&\\
b_i=&0&1&2&3&   &\{b_i\}=5\\
b_i'=&0&2&4&6&(\lambda=2;\delta=0)&\\
c_i=&0&0&2&6&   &\{c_i\}=6\\
c_i'=&4&4&8&16&(\lambda=2;\delta=4)&\\
d_i=&0&4&8&16&   &\{d_i\}=7\\
&&&&etc.&&
\end{matrix}$$
In actual use, it is not necessary to build up one step at a time.  Several shortcuts may be devised.  One of them, suggested by Prof. J. Ginsburg, enables us to skip upward in a doubling fashion.  In the formulae for $b_i'$ above, eliminate the $a_i$ using $(17)$ with $\lambda=2, \delta=0$.  Then
$$\begin{matrix}
b_2=&-a_2&-a_3&+a_4\\
b_3=&&-2a_3&+2a_4\\
b_4=&-a_2&-a_3&+3a_4.
\end{matrix}$$

At this point it will be helpful to change our symbolism, putting
$$\begin{matrix}
A=a_2&  &A_1=b_2\\
B=a_3&  &B_1=b_3\\
C=a_4 && C_1=b_4.
\end{matrix}$$
Then the set $(0,A_1,B_1,C_1)$ lasts one step longer than $(0,A,B,C)$ and they are connected by the equations
$$\begin{matrix}
A_1=&-A&-B&+&C\\
B_1=&&-2B&+&2C\\
C_1=&-A&-B&+&3C.
\end{matrix}$$
Applying the same transformation to $A_1, B_1, C_1$, we get
$$A_2=-A_1-B_1+C_1=2B$$
$$B_2=-2B_1+2C_1=-2A+2B+2C$$
$$C_2=-A_1-B_1+3C_1=-2A+6C.$$
The factor $2$ can be removed without affecting the \lq\lq life" of the set.  We write
$$\begin{matrix}
A_2=B&&\\
B_2=-A&&+B&+C\\
C_2=-A&&+3C.
\end{matrix}$$
Similarly,
$$\begin{matrix}
A_4=B_2& & &=-A+B+C\\
B_4=-A_2&+B_2&+C_2&=-2A+2B+4C\\
C_4=-A_2&+3C_2&&=-3A-B+9C.
\end{matrix}$$
This last transformation increases the \lq\lq life" of a set by four steps.  It is clear that by repetition we can get formulae for building up $8,16,32,$ etc., at a time.  Using the formulae for a skip upward of $16$, and setting $A=0, B=0, C=1$, we get the set $(0,193,548,1201)$ which lasts $16$ more sets than $(0,0,0,1)$ and therefore lasts a total of $20$ differencings.

Another approach to the problem of building up several steps at a time is to consider the formulae for building up one stage as difference equations:
$$\begin{matrix}
(1)&A_{n+1}=&-A_n &-B_n&+C_n\\
(2)& B_{n+1}=&&-2B_n&+2C_n\\
(3)&C_{n+1}&-A_n&-B_n&+3C_n.
\end{matrix}$$

These may be simplified as follows,
$$\begin{matrix}
(4)& C_{n+1}-A_{n+1}=2C_n&  (3)-(1)&\\
(5)& B_{n+1}-2A_{n+1}=2A_n & (2)-(1)-(1)&\\
(6)& A_{n+1}=C_{n+1}-2C_n & (4)&
\end{matrix}$$
$$(7)\; B_{n+1}=2A_{n+1}+2A_n  =2C_{n+1}-4C_n+2C_n-4C_{n-1}\; (5) \; \textrm{and}\; (6)$$
$$(8)\; C_{n+1}=-C_{n}+2C_{n-1}-2C_n+2C_{n-1}+4C_{n-2}+3C_n,\;\; \; \textrm{putting\; (6)\; and \;(7)\; in\; (3)}$$
$$\left.\begin{matrix}
(9)& C_{n+1}=4C_{n-1}+4C_{n-2}. & & \\
\textrm{Similarly} &&&\\
(10)& B_{n+1}=4B_{n-1}+4B_{n-2} & &\\
(11)& A_{n+1}=4A_{n-1}+4A_{n-2}. & &
\end{matrix}\right\}$$

Given $A_0\; B_0\; C_0, A_1, B_1, C_1$ may be calculated by $(1), (2), (3), $ and $A_2\; B_2\; C_2,$ etc., worked out by successive applications of $(9), (10), (11)$.

To generalize the last result we must consider what the symbol $\{a_i\}$ means if $k$ is not a power of $2$.  Let us stipulate then that the game terminates just short of the first repetition.  E.g., if $A_0=2,5,9$ as in the example attached to Theorem II, $\{a_0\}$ will be $7$, as $A_8$ repeats $A_5$.  With this in mind we can investigate the possibility of building upward for any $k$ as we have just done for $k=4$.

First observe that equations $(12), (13),$ and $(14)$ are valid for any $k$ under our extended definition.  $(16)$ is equally general.

We will exam three cases:

\begin{itemize}
\item[(I)] $k=2$
\item[(II)] $k=$ an odd number $>2$
\item[(III)] $k=$ an even number $>2$.
\end{itemize}

Case $(I)$: If $k=2$,
$$A_0=a_1 \;\;\; a_2$$
$$A_1=|a_2-a_1| |a_1-a_2|$$
$$A_2=0\;\;\;0.$$
So that $\{a_i\}\leq 2 \;\;\; i=1,2.$ I.e., the building upward process is impossible for $k=2$.

Case $II$: $ k=2l+1, l>0$.

Take $S=0,a_1,b_1,a_2,b_2,...,a_l,b_l.$  Clearly there are $2l+1$ or $k$ members.  Suppose that $a_i$ and $b_i$ satisfy the following conditions:
$$(\alpha')\;\;\;  a_i\; \textrm{and}\; b_i \; \textrm{are postive integers}$$
$$ (\beta')\;\;\; b_i> a_i\;\;\; i=1,2,... l$$
$$(\gamma')\;\;\; b_l>b_i$$
$$\;\;\;  a_l> a_i \;\;\;i=2,3,...l$$
Set
$$\left\{\begin{matrix}
\delta& =&\sum_{i=1}(b_i-a_i)& (\textrm{note that}\; \delta>0)\\
A_1& =&\delta& \\
A_i& =&\delta -\sum_{j=1}^{i-1}(b_j-a_j)& i=2,3,...l\\
B_1& =&2\delta+a_1 & \\
B_i&=& 2\delta+a_i-\sum_{j=1}^{i-1}(b_j-a_j) & i=2,3,...l\\
T& =&0,&A_1, B_1, A_2,B_2, & ...,A_l, B_l
\end{matrix}\right.$$
We shall demonstrate that $\{T\}=\{S \}+1.$  First we prove that $|\Delta|T=S+\delta.$
$$|\Delta|0=|A_1-0|=A_1=\delta$$
$$|\Delta|A_i=|B_i-A_i|=\delta+a_i \;\;\; i=1,2,3,...l.$$
$$|\Delta|B_i=|A_{i+1}-B_i|=\delta+b_i \;\;\; i=1,2,3,...(l-1)$$
$$|\Delta|B_l=|0-B_l|=B_l=2\delta+a_i-\sum_{1}^{l-1}(b_j-a_j)=\delta+b_l (\; \textrm{since} \;\delta=\sum_{1}^l(b_j-a_j))$$
$$\therefore |\Delta|T=S+\delta.$$

Now, $B_1=2\delta+a_1=\delta+a_1+\sum_{1}^l(b_j-a_j)=\delta+b_1+\sum_{2}^l(b_j-a_j).$

By virtue of condition $\beta', B_1>b_1+\delta$.  By condition $\gamma'$, this makes $B_1$ greater than every member of $(S+\delta)$.  Therefore the set $T$ cannot be repeated in the derived sets of $(S+\delta)$, and we are justified in saying that
\begin{equation}
\{T\}=\{S+\delta\}+1
\end{equation}
Or, using $(13)$
\begin{equation}
\{T\}=\{S\}+1
\end{equation}
All that remains is to verify that every new set $T$ satisfies the restrictions $\alpha', \beta', \gamma'$.

For $(\alpha')$:
$$A_1=\delta>0$$
$$A_i=\delta-\sum_{j=1}^{i-1}(b_j-a_j)=\sum_i^l(b_j-a_j)>0 \;\;\; i=2,3,...l$$
$$B_1=2\delta+a_1>0$$
$$B_i=2\delta+a_l-\sum_{j=1}^{i-1}(b_j-a_j)=\delta+a_i+\sum_{i}^l(b_j-a_j)>0 \;\;\; i=2,3,...l.$$

For $(\beta')$:
$$B_i-A_i=a_i+\delta>0 \;\;\; i=1,2,3,...l.$$
For $(\gamma')$:
$$B_1-B_i=a_l-a_i+\sum_{j=1}^{i-1}(b_j-a_j)>0\;\;\; i=2,3,...l$$
$$A_1-A_i=\sum_{j=1}^{i-1}(b_j-a_j)>0 i=2,3,...l.$$

As $T$ satisfies the same conditions we imposed on $S$, a new set can be built up on $T$ in the same way.  By repeating the process, it is possible to construct a set $U$ with an odd number of members, for which $\{U\}$ is as large as desired.

Case $III$: $k=2l+2, l>0.$

Take $S=0,c,a_a,b_1,a_2,b_2,...a_l,b_l.$  Suppose that $a_i$ and $b_i$ satisfy
$$ (\alpha'')\;\;\; a_i\geq 0; b_i\geq 0; c\geq 0$$
$$(\beta'')\;\;\; b_1> a_1+c; b_i>a_i \;\;\;i=2,3,...,l$$
$$(\gamma'')\;\;\; b_1>b_i; a_1>a_i\;\;\; i=2,3,...,l.$$
Set
$$\left\{\begin{matrix}
C&=&-c+\sum_{j=1}^l(b_j-a_j)& \\
A_i& =&2[\sum_{j=i}^l(b_j-a_j)]& i=1,2,,3,...l\\
B_i&= &2[\sum_{j=i}^l(b_j-a_j)]+2a_i-c+\sum_{j=1}^l(b_j-a_j)& i=2,3,..l\\
T& =&0,C, A_1, B_1, A_2,B_2,...A_l, B_l.&
\end{matrix}\right.$$

In a parallel manner to the proof of Case $\textrm{II}$, we demonstrate that
$$|\Delta|T=[2S-c+\sum_{j=1}^l(b_j-a_j)]$$
and that
\begin{equation}
\{T\}=\{S\}+1\
\end{equation}
Also $T$ satisfies the restrictions $\alpha'', \beta'', \gamma''$, and the reasoning of Case $\textrm{II}$ applies here.

We can combine the results of Case $\textrm{I}$, $\textrm{II}$, and $\textrm{III}$ into
\begin{theorem}
If $k$ is the number of members of the original set, then for $k>2$, there is no upper limit to $\{S_0\}$.
\end{theorem}

To summarize, in the $k$-number game, $k$ represents the number of numbers of the original set, $\{S_0\}$ the number of differencings before a set repeats.  Then if $k$ is a power of $2$, $S_{\{S_0\}}$ will be composed of all zeros.  If $k$ is not a power of $2$, this will not happen in general.  And if $k$ is greater than $2$, sets may be found for any value of $\{S\}$ no matter how high.

\end{document}